\documentclass[11pt,twoside]{article}

\usepackage{graphicx}
\usepackage{a4wide}
\usepackage{amsfonts}
\usepackage{amssymb}
\usepackage{amsmath}
\usepackage{graphicx}

\newcommand{\ignore}[1]{}

\def\@begintheorem#1#2{\par\bgroup{\sc #1\ #2. }\it\ignorespaces}
\def\@opargbegintheorem#1#2#3{\par\bgroup{\sc #1\ #2\ (#3). } \it\ignorespaces}
\def\@endtheorem{\egroup}
\newtheorem{theorem}{Theorem}[section]
\newtheorem{corollary}[theorem]{Corollary}
\newtheorem{lemma}[theorem]{Lemma}

\newtheorem{example}[theorem]{Example}
\newtheorem{proposition}[theorem]{Proposition}
\newtheorem{definition}[theorem]{Definition}
\newcommand{\bt}[1]{\begin{theorem}\label{#1}}
\newcommand{\bc}[1]{\begin{corollary}\label{#1}}
\newcommand{\bl}[1]{\begin{lemma}\label{#1}}
\newcommand{\be}[1]{\begin{example}\label{#1}}
\newcommand{\bp}[1]{\begin{proposition}\label{#1}}
\newcommand{\ba}[1]{\begin{algorithm}\rm\label{#1}}
\newcommand{\bd}[1]{\begin{definition}\rm\label{#1}}{\normalsize }
\newcommand{\bpr}{\noindent {\em Proof. }}
\newcommand{\et}{\end{theorem}}
\newcommand{\ec}{\end{corollary}}
\newcommand{\el}{\end{lemma}}
\newcommand{\ee}{\end{example}}
\newcommand{\ep}{\end{proposition}}
\newcommand{\ed}{\end{definition}}
\newcommand{\epr}{{\ \vbox{\hrule\hbox{%
\vrule height1.3ex\hskip0.8ex\vrule}\hrule}}\\\par}

\def\Z{\mathbb{Z}}

\begin{document}

\title{\bf Uniform and Monotone Line Sum Optimization}

\author{
Martin Kouteck\'y
\thanks{\small Computer Science Institute, Faculty of Mathematics and Physics,
Charles University. Email: koutecky@iuuk.mff.cuni.cz}
\and
Shmuel Onn
\thanks{\small Technion - Israel Institute of Technology.
Email: onn@technion.ac.il}
}
\date{}

\maketitle

\begin{abstract}
The {\em line sum optimization problem} asks for a $(0,1)$-matrix
minimizing the sum of given functions evaluated at its row and column sums.
We show that the {\em uniform} problem, with identical row functions and
identical column functions, and the {\em monotone} problem, over matrices
with nonincreasing row and column sums, are polynomial time solvable.

\vskip.2cm
\noindent {\bf Keywords:}
majorization, column sum, row sum, matrix, degree sequence, graph

\vskip.2cm
\noindent {\bf MSC:}
05A, 15A, 51M, 52A, 52B, 52C, 62H, 68Q, 68R, 68U, 68W, 90B, 90C
\end{abstract}

\section{Introduction}

For a positive integer $n$ let $[n]:=\{1,2,\dots,n\}$. For an $m\times n$ matrix $A$ let
$r_i(A):=\sum_{j=1}^n A_{i,j}$ for $i\in[m]$ be its row sums
and let $c_j(A):=\sum_{i=1}^m A_{i,j}$ for $j\in[n]$ be its column sums.

\vskip.2cm
We consider here the following algorithmic problem.

\vskip.2cm\noindent{\bf Line Sum Optimization.} Given $m,n$ and functions
$f_i:\{0,1,\dots,n\}\rightarrow\Z$ for $i\in[m]$ and $g_j:\{0,1,\dots,m\}\rightarrow\Z$
for $j\in[n]$, find an $m\times n$ $(0,1)$-matrix, $A\in\{0,1\}^{m\times n}$, which minimizes
$$\sum_{i=1}^m f_i(r_i(A))+\sum_{j=1}^n g_j(c_j(A))\ .$$

\vskip.2cm
For instance, for $m=n=4$, and functions $f_i(x)=(x-1)^2(x-3)^2$ for $i\in[4]$ and
$g_j(x)=(x-2)^2(x-3)^2$ for $j\in[4]$, an optimal solution is the following matrix,
with row sums $(r_1,r_2,r_3,r_4)=(3,3,3,1)$ and column sums
$(c_1,c_2,c_3,c_4)=(3,3,2,2)$ and objective value $0$,
$$A\ = \left(\begin{array}{cccc}
1&1&1&0\\
1&1&1&0\\
1&1&0&1\\
0&0&0&1\\
\end{array}
\right)\ .$$

\vskip.2cm
In this article we consider the {\em uniform} case, where all $f_i$ are the same
and all $g_j$ are the same, and show that the problem can then be solved in polynomial time.
\bt{uniform}
The uniform line sum optimization problem, where for some given functions $f,g$ we have
$f_i=f$ for all $i\in[m]$ and $g_j=g$ for all $j\in[n]$, can be solved in polynomial time.
\et

We call a matrix $A$ {\em monotone} if it has nonincreasing row and column sums,
that is, $r_1\geq\cdots\geq r_m$ and $c_1\geq\cdots\geq c_n$.
We also solve the line sum problem over monotone matrices.

\bt{monotone}
Given $m,n$, $f_i:\{0,1,\dots,n\}\rightarrow\Z$, and $g_j:\{0,1,\dots,m\}\rightarrow\Z$,
a monotone $A\in\{0,1\}^{m\times n}$ minimizing
$\sum_{i=1}^m f_i(r_i(A))+\sum_{j=1}^n g_j(c_j(A))$ is polynomial time computable.
\et

Theorem \ref{monotone} clearly implies Theorem \ref{uniform}: if for some given
functions $f,g$ we have $f_i=f$ for all $i\in[m]$ and $g_j=g$ for all $j\in[n]$,
then the objective value of any matrix is invariant under row and column permutations,
and hence an optimal solution to the monotone problem is an optimal solution
to the uniform problem as well. So we need only prove Theorem \ref{monotone}.

\vskip.2cm
The {\em uniform column sum problem}, where the row sums $r_1,\dots,r_m$ are specified, and the
objective is to minimize $\sum_{j=1}^n g(c_j(A))$, recently solved in \cite{Onn}, is a special
case of Theorem \ref{monotone}, obtained by assuming $r_1\geq\cdots\geq r_m$ and taking
$f_i(x)=a(x-r_i)^2$ for all $i$ and sufficiently large $a$. The line sum problem is a special
case of the {\em degree sequence optimization problem}, where, given a graph $H=(V,E)$ and
functions $f_v:\{0,1,\dots,d_v(H)\}\rightarrow\Z$ for $v\in V$, with $d_v(H)$ the degree of $v$
in $H$, we need to find a subgraph $G=(V,F)\subseteq H$ minimizing $\sum_{v\in V}f_v(d_v(G))$.
Indeed, identifying matrices $A\in\{0,1\}^{m\times n}$ with bipartite graphs $G=(V,F)$ where
$V=\{u_1,\dots,u_m\}\uplus\{w_1,\dots,w_n\}$ and $F=\{\{u_i,w_j\}\,:\,A_{i,j}=1\}$,
the line sum problem reduces to the degree sequence problem with $H=K_{m,n}$ the complete
bipartite graph. In the case of $H=K_n$ the complete graph, the {\em uniform problem},
where all functions are the same, $f_v=f$ for all $v\in V$, was recently shown in
\cite{DLMO} to be polynomial time solvable, using the characterization of degree sequences
by Erd\H{o}s and Gallai \cite{EG}. For general graphs $H$, the problem was shown in \cite{AS}
to be NP-hard already when $f_v(x)=-x^2$ for all $v\in V$, but is polynomial time solvable
if the functions are convex \cite{AS,DO}. We conjecture that the degree sequence problem over
$H=K_{m,n}$, which is the line sum problem, as well as over $H=K_n$, is polynomial time
solvable for arbitrary functions at the vertices, not necessarily identical.

\section{Proof}
\label{proof}

For an $m\times n$ matrix $A$ let $r(A)=(r_1(A),\dots,r_m(A))$ and $c(A)=(c_1(A),\dots,c_n(A))$
be the tuples of row and column sums, and let $f(r(A))=\sum_{i=1}^m f_i(r_i(A))$ and
$g(c(A))=\sum_{j=1}^n g_j(c_j(A))$. We need the following terminology. A nonincreasing
$r=(r_1,\dots,r_m)$ is {\em majorized} by a nonincreasing $s=(s_1,\dots,s_m)$ if
$\sum_{i=1}^h r_i\leq\sum_{i=1}^h s_i$ for $h\in[m]$ and $\sum_{i=1}^m r_i=\sum_{i=1}^m s_i$.
(See \cite{MOA} for more details on the theory and applications of majorization.)
The {\em conjugate} of a nonincreasing tuple $c=(c_1,\dots,c_n)$ with $c_1\leq m$ is
the nonincreasing tuple $s=(s_1,\dots,s_m)$ where $s_i=|\{j\,:\,c_j\geq i\}|$ for $i\in[m]$.
Note that $s_1\leq n$ and $\sum_{i=1}^m s_i=\sum_{j=1}^n c_j$.
We make use of the following characterization due to Ryser \cite{Rys}.

\bp{Ryser}
A monotone $A\in\{0,1\}^{m\times n}$ with row and column sums
$r=(r_1,\dots,r_m)$ and $c=(c_1,\dots,c_n)$ exists if and only
if $r$ is majorized by the conjugate $s=(s_1,\dots,s_m)$ of $c$.
\ep

\begin{figure}
	\begin{center}
\includegraphics[width=0.7\textwidth]{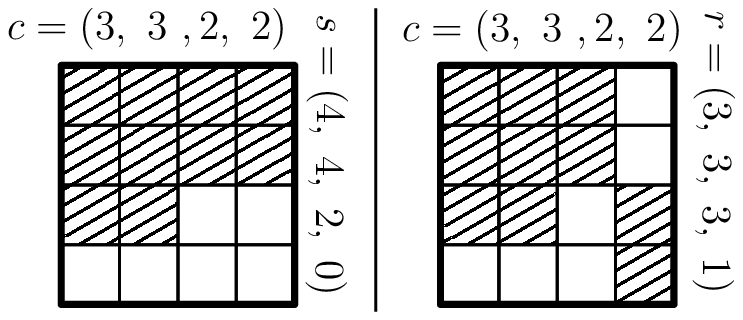}
\end{center}
\caption{A nonincreasing tuple $c = (3,3,2,2)$, its conjugate $s=(4,4,2,0)$, a nonincreasing tuple $r=(3,3,3,1)$, which is majorized by $s$, and a monotone matrix $A\in\{0,1\}^{4\times 4}$ with columns sums $c$ and row sums $r$.} \label{fig:1}
\end{figure}

For instance, if $m=n=4$ and $c=(3,3,2,2)$, then $s=(4,4,2,0)$, and $r=(3,3,3,1)$ is majorized
by $s$ so there is a monotone matrix $A\in\{0,1\}^{4\times 4}$ with row sums $r$ and column sums $c$,
see Figure~\ref{fig:1}.

\vskip.2cm
Next we note that a matrix with given row and column sums can be efficiently obtained.
\bl{flow}
Given $r=(r_1,\dots,r_m)$ and $c=(c_1,\dots,c_n)$ which satisfy the Ryser criterion, a matrix
$A$ in $\{0,1\}^{m\times n}$ with row sums $r$ and column sums $c$ is polynomial time computable.
\el
\bpr
The problem is solvable either by the efficient simple Gale-Ryser algorithm, see
\cite[Chapter 3]{Bru}, or by network flows as follows. Define a directed graph with capacities
on the edges as follows. There are vertices $s,t$, $u_1,\dots,u_m$, and $w_1,\dots,w_n$.
There are edges $[s,u_i]$ for $i\in[m]$ with capacity $r_i$, edges $[u_i,w_j]$ for $i\in[m]$
and $j\in[n]$ with capacity $1$, and edges $[w_j,t]$ for $j\in[n]$ with capacity $c_j$.
Then, as is well known, a maximum nonnegative integer flow from $s$ to $t$ can be computed
in polynomial time, see e.g. \cite{Sch}. Then $A$ is read off from the maximum flow by
taking $A_{i,j}$ to be the flow on edge $[u_i,w_j]$ for all $i\in[m]$ and $j\in[n]$.
\epr

A key idea facilitating our algorithm is a different way to view or encode the nonincreasing tuples $c=(c_1,\dots,c_n)$ with $c_1\leq m$ and their conjugates $s=(s_1,\dots,s_m)$ with $s_1\leq n$.
Viewing a nonincreasing tuple $c$ as a series of left-aligned and top-aligned ``strips'' (see Figure~\ref{fig:1}).
In this view, $c$ describes the length of each strip.
An alternative description is by specifying the strip lengths which occur together with their multiplicities.
As we will see, this view is much more amenable to designing a dynamic program.
Let us now describe it in detail.

\begin{figure}
	\begin{center}
		\includegraphics[width=0.7\textwidth]{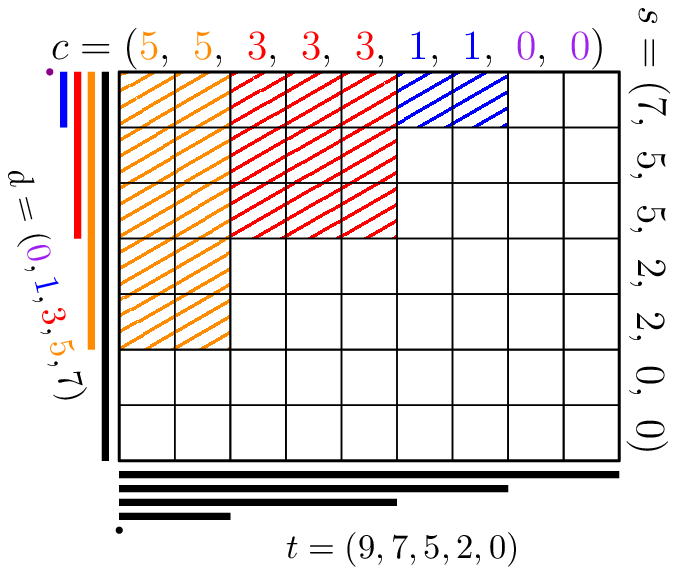}
	\end{center}
	\caption{The column sums tuple $c$ and its conjugate $s$
constructed in the path of Example~\ref{example_bipartite}.} \label{fig:3}
\end{figure}

The {\em type} of a nonincreasing tuple is the number of
distinct nonzero values among its components, i.e., the number of occurring strip lengths.
It is easy to see that if $s$ is the
conjugate of $c$ then $s$ and $c$ have the same type $0\leq k\leq\min\{m,n\}$.
For a tuple $c$ of type $k$, we define numbers $n=t_0\geq t_1>\cdots>t_k>t_{k+1}=0$ and $0=d_0<d_1<\cdots<d_k\leq d_{k+1}=m$ such that, for $h=0,1,\dots,k$, we have that $c$ has $t_h - t_{h+1}$ components equal to $d_h$.
Clearly each $c$ has such an encoding using vectors $t=(t_0, \dots, t_{k+1})$ and $d=(d_0, \dots, d_{k+1})$, and also each choice of vectors $t$ and $d$ as previously mentioned corresponds to exactly one $c$.
Moreover, the conjugate $s$ of $c$ encoded using $t$ and $d$ has $d_{h+1}-d_h$ components equal to $t_{h+1}$
for $h=0,1,\dots,k$. A particularly neat way to describe $c$, $s$ using $t$, $d$ is the abridged form:
\begin{equation}\label{type}
c=(c_1,\dots,c_n)=(d_k^{t_k-t_{k+1}},\dots,d_0^{t_0-t_1})\ ,\quad
s=(s_1,\dots,s_m)=(t_1^{d_1-d_0},\dots,t_{k+1}^{d_{k+1}-d_k})\enspace .
\end{equation}
Continuing with the example before we have $c=(3,3,2,2) = (3^2, 2^2)$ and $s = (4,4,2,0) = (4^2,2^1,0^1)$.
See also Figure \ref{fig:3} for a larger example which will be treated in detail in Example \ref{example_bipartite}. 
Note that any $c$ and its conjugate $s$ of type $k$ arise that way from some such $t$, $d$.

\vskip.2cm
We can now solve the monotone and hence also the uniform line sum problems.

\vskip.2cm\noindent
{\bf Theorem \ref{monotone}}
{\em  Given $m,n$, $f_i:\{0,1,\dots,n\}\rightarrow\Z$, and
$g_j:\{0,1,\dots,m\}\rightarrow\Z$, a monotone $A\in\{0,1\}^{m\times n}$ minimizing
$\sum_{i=1}^m f_i(r_i(A))+\sum_{j=1}^n g_j(c_j(A))$ is polynomial time computable.}

\vskip.2cm
\bpr
By Proposition \ref{Ryser} and Lemma \ref{flow} it suffices to find nonincreasing tuples\break
$c=(c_1,\dots,c_n)$ with $c_1\leq m$ and $r=(r_1,\dots,r_m)$ majorized by the conjugate
$s=(s_1,\dots,s_m)$ of $c$ that minimize $f(r)+g(c)=\sum_{i=1}^m f_i(r_i)+\sum_{j=1}^n g_j(c_j)$.
For this, we use the encoding of tuples $c$ of type $k$ and their conjugates $s$ discussed
above. For type $k=0$ we trivially have $c=(0^n)$ and $r=s=(0^m)$ with value
$\sum f_i(0)+\sum g_j(0)$. For each type $k=1,\dots,\min\{m,n\}$ we provide a construction which
reduces the problem of finding the best $c,r$ where $c$ has type $k$ to that of finding a
shortest directed path in a suitable directed graph $D_k$ with lengths on the edges.
An alternative perspective is that this is a dynamic programming algorithm where we gradually solve larger and larger subproblems; we choose the shortest path encoding to focus attention on the information we keep in each state in order to compute the next one.

\vskip.2cm
We now describe $D_k$. There are two special vertices $u,v$. The remaining vertices  are labeled by septuples
of integers $(h,t_h,d_h,i,r_i,S_i,R_i)$, where $h\in[k+1]$, $0\leq t_h,r_i\leq n$,
$d_h,i\in[m]$, $S_i,R_i\in[mn]$. We always define the ``boundary" values
$t_0=n$, $d_{k+1}=m$, and $d_0=t_{k+1}=0$. Our goal is to encode each column vector $c$ of type $k$, implicitly its conjugate $s$, and each $r$ majorized by $c$,
in a directed $u-v$ path of length $f(r)+g(c)$.

\vskip.2cm
Before formally describing the edges and their lengths, we explain how such $c,r$ give a path.
The reader is referred to Example \ref{example_bipartite} below for a specific demonstration.
Consider any choice of numbers $n=t_0\geq t_1>\cdots>t_k>t_{k+1}=0$ and
$0=d_0<d_1<\cdots<d_k\leq d_{k+1}=m$. These numbers define the tuple $c=(c_1,\dots,c_n)$ and
its conjugate $s=(s_1,\dots,s_m)$ as in \eqref{type}, where $c$ has $t_h-t_{h+1}$ components
equal to $d_h$ and $s$ has $d_{h+1}-d_h$ components equal to $t_{h+1}$ for $h=0,1,\dots,k$.
Now consider any choice of a nonincreasing tuple $r=(r_1,\dots,r_m)$ majorized by $s$.
Let $S_i=\sum_{j=1}^i s_i$ and $R_i=\sum_{j=1}^i r_i$. For $r$ to be nonincreasing we need
$r_{i+1}\leq r_i$ for $1\leq i<m$. For $r$ to be majorized by $s$ we need $R_i\leq S_i$ for
$i\in[m]$ and $R_m=S_m$. The path corresponding to such choices
(see Example \ref{example_bipartite}) is as follows. It has $m+2$ vertices, starting at $u$,
going through $m$ vertices $(h,t_h,d_h,i,r_i,S_i,R_i)$ with $i=1,\dots,m$ and
$r_1,\dots,r_m$ the components of the chosen tuple $r$, and ending at $v$.

\begin{figure}
	\begin{center}
		\includegraphics[width=0.7\textwidth]{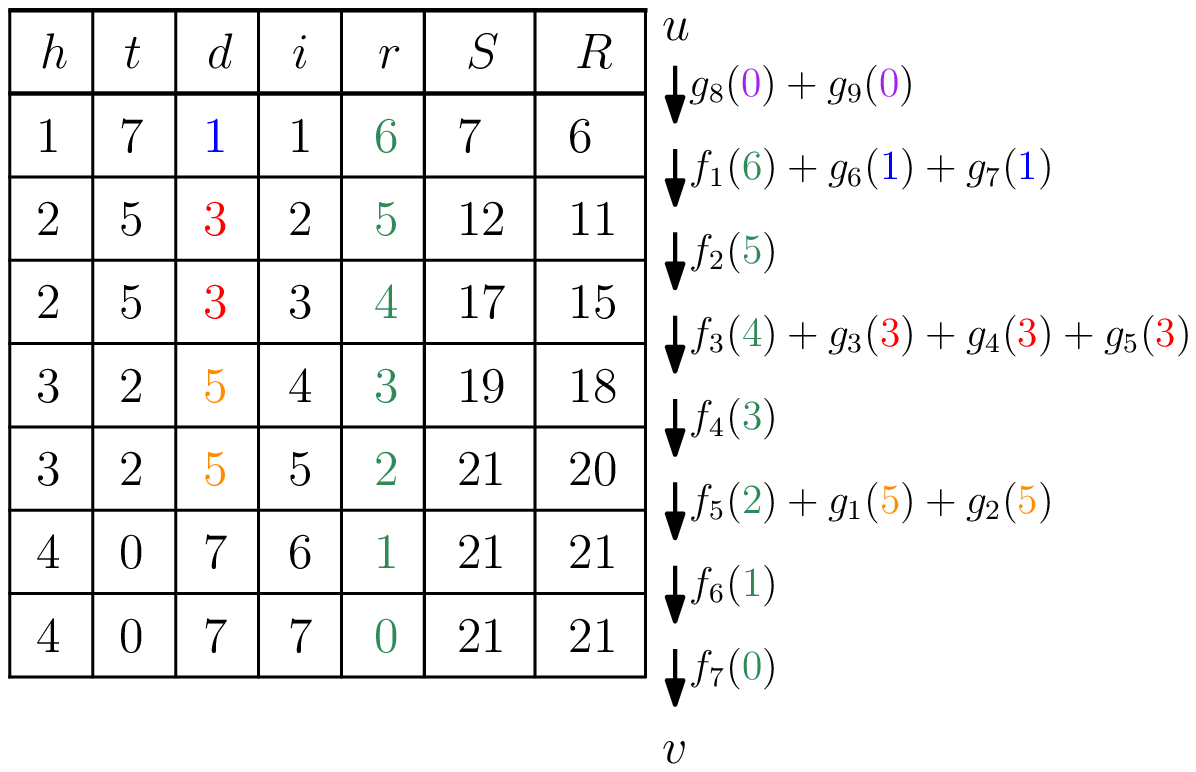}
	\end{center}
	\caption{The path of Example~\ref{example_bipartite}. Each row is one internal vertex of the path; the values on the right of the table are the lengths
of the edges between the corresponding consecutive vertices. Notice that $t$, $d$ are such that they encode some nonincreasing $c$; $S$ are prefix sums of the conjugate $s$ of $c$, and because in each row we have
$S \geq R$, $R$ are prefix sums of some nondecreasing $r$ which is majorized by $s$.} \label{fig:2}
\end{figure}

More specifically, we start with vertex $u$ and go to vertex
$$(h=1,t_1,d_1,i=d_0+1=1,r_1\leq t_1,S_1=t_1,R_1=r_1\leq S_1)$$
along an edge of length $\sum\{g_j(0):t_1<j\leq t_0\}$ accounting for the contribution
of the $(t_0-t_1)=n-t_1$ components $d_0=0$ of $c$ if any. We proceed on a
path where $h=1,t_1,d_1$ remain fixed, while we increment $i$ from $d_0+1=1$
to $d_1$, where the components $s_{d_0+1},\dots,s_{d_1}$ of $s$ are all equal to $t_1$
so that we set $S_{i+1}=S_i+t_1$ for their sum. The components $r_{d_0+1},\dots,r_{d_1}$
of $r$ are as chosen and we set $R_{i+1}=R_i+r_{i+1}$ for their sum. If $d_1>d_0+1$ then
the length of the edge from the vertex with index $i\geq d_0+1$ to $i+1\leq d_1$ is
$f_i(r_i)$ accounting for the contribution of $r_i$.
When $i$ reaches $d_1$, we increment $h$ and proceed to vertex
$$(h=2,t_2,d_2,i+1=d_1+1,r_{i+1}\leq r_i,S_{i+1}=S_i+t_2,R_{i+1}=R_i+r_{i+1}\leq S_{i+1})$$
along an edge of length $f_{d_1}(r_{d_1})+\sum\{g_j(d_1):t_2<j\leq t_1\}$ accounting for
the contribution of $r_{d_1}$ and the $(t_1-t_2)$ components $d_1$ of $c$. Now we fix
$h=2,t_2,d_2$ and continue on a path where we increment $i$ from $d_1+1$ to $d_2$, where
the components $s_{d_1+1},\dots,s_{d_2}$ of $s$ are all equal to $t_2$ so that we set
$S_{i+1}=S_i+t_2$ for their sum. The components $r_{d_1+1},\dots,r_{d_2}$ of $r$ are as
chosen and we set $R_{i+1}=R_i+r_{i+1}$ for their sum. We continue this way till we arrive
at the vertex $(h=k,t_k,d_k,i=d_k,r_{d_k},S_{d_k},R_{d_k}\leq S_{d_k})$.
If $i=d_k=m$ and $R_m=S_m$ then we move to $v$ along an edge of length
$f_m(r_m)+\sum\{g_j(m):t_{k+1}<j\leq t_k\}$ accounting for the contribution of $r_m$
and the $(t_k-t_{k+1})=t_k$ components $d_k=m$ of $c$. If $i=d_k<m$ then we move to
$$(h=k+1,t_{k+1}=0,d_{k+1}=m,i+1=d_k+1,r_{i+1}\leq
r_i,S_{i+1}=S_i,R_{i+1}=R_i+r_{i+1}\leq S_{i+1})$$
along an edge of length $f_{d_k}(r_{d_k})$ accounting for the contribution of $r_{d_k}$.
We proceed on a path where $h=k+1,t_{k+1}=0,d_{k+1}=m$ remain fixed, while we increment
$i$ from $d_k+1$ to $m$, where the components $s_{d_k+1},\dots,s_m$ of $s$ are all equal
to $t_{k+1}=0$ so that $S_{i+1}=S_i$. The components $r_{d_k+1},\dots,r_m$ of $r$ are as
chosen and we set $R_{i+1}=R_i+r_{i+1}$ for their sum. If $m=d_{k+1}>d_k+1$ then the length of
the edge from the vertex with index $i\geq d_k+1$ to $i+1\leq m$ is $f_i(r_i)$ accounting for
the contribution of $r_i$. Finally, we arrive at the vertex $(h=k+1,0,m,i=m,r_m,S_m,R_m)$,
and if $R_m=S_m$ then we move to $v$ along an edge of length $f_m(r_m)$ accounting
for the contribution of $r_m$. Let us now work through an example with $d_k<m$.
After the example we will complete the formal description of $D_k$.

\be{example_bipartite}
We now demonstrate the construction of the directed graph $D_k$. Consult also Figures \ref{fig:2} and \ref{fig:3}.
Let $m=7$, $n=9$,
$c=(5,5,3,3,3,1,1,0,0)$, and $r=(6,5,4,3,2,1,0)$. The conjugate of $c$ is $s=(7,5,5,2,2,0,0)$ which
majorizes $r$. The type of $c$ and $s$ is $k=3$. The tuples $t=(t_0,\dots,t_{k+1})=(9,7,5,2,0)$
and $d=(d_0,\dots,d_{k+1})=(0,1,3,5,7)$ define $(5^2,3^3,1^2,0^2)=c$ and $(7,5^2,2^2,0^2)=s$
as in \eqref{type}. The directed $u-v$ path in $D_3$ corresponding
to $c$ and $r$, with edge lengths indicated (see also Figure \ref{fig:2}), is:
$$u$$
$$\downarrow \sum\{g_j(d_0):t_1<j\leq t_0\}=g_8(0)+g_9(0)$$
$$(h=1,t_1=7,d_1=1,i=1=d_1,r_1=6\leq t_1,S_1=t_1=7,R_1=r_1=6\leq S_1)$$
$$\downarrow f_1(r_1)+\sum\{g_j(d_1):t_2<j\leq t_1\}=f_1(6)+g_6(1)+g_7(1)$$
$$(h=2,t_2=5,d_2=3,i=2<d_2,r_2=5\leq r_1,S_2=S_1+t_2=12,R_2=R_1+r_2=11\leq S_2)$$
$$\downarrow f_2(r_2)=f_2(5)$$
$$(h=2,t_2=5,d_2=3,i=3=d_2,r_3=4\leq r_2,S_3=S_2+t_2=17,R_3=R_2+r_3=15\leq S_3)$$
$$\downarrow f_3(r_3)+\sum\{g_j(d_2):t_3<j\leq t_2\}=f_3(4)+g_3(3)+g_4(3)+g_5(3)$$
$$(h=3,t_3=2,d_3=5,i=4<d_3,r_4=3\leq r_3,S_4=S_3+t_3=19,R_4=R_3+r_4=18\leq S_4)$$
$$\downarrow f_4(r_4)=f_4(3)$$
$$(h=3,t_3=2,d_3=5,i=5=d_3,r_5=2\leq r_4,S_5=S_4+t_3=21,R_5=R_4+r_5=20\leq S_5)$$
$$\downarrow f_5(r_5)+\sum\{g_j(d_3):t_4<j\leq t_3\}=f_5(2)+g_1(5)+g_2(5)$$
$$(h=4,t_4=0,d_4=7,i=6<d_4,r_6=1\leq r_5,S_6=S_5+t_4=21,R_6=R_5+r_6=21\leq S_5)$$
$$\downarrow f_6(r_6)=f_6(1)$$
$$(h=4,t_4=0,d_4=7,i=7=d_4,r_7=0\leq r_6,S_7=S_6+t_4=21,R_7=R_6+r_6=21=S_5)$$
$$\downarrow f_7(r_7)=f_7(0)$$
$$v$$
So the total length of this path is indeed equal to the
objective value corresponding to $r$ and $c$,
\begin{eqnarray*}
	&&\left(f_1(6)+f_2(5)+f_3(4)+f_4(3)+f_5(2)+f_6(1)+f_7(0)\right)\\
	&+&\left(g_1(5)+g_2(5)+g_3(3)+g_4(3)+g_5(3)+g_6(1)+g_7(1)+g_8(0)+g_9(0)\right)\ =\ f(r)+g(c)\ .
\end{eqnarray*}
Note that by Proposition \ref{Ryser} a matrix with sums $r,c$ exists and can be found
by Lemma \ref{flow},
$$A\ = \left(\begin{array}{ccccccccc}
	1&1&1&1&1&1&0&0&0\\
	1&1&1&1&1&0&0&0&0\\
	1&1&1&1&0&0&0&0&0\\
	1&1&0&0&0&0&1&0&0\\
	1&1&0&0&0&0&0&0&0\\
	0&0&0&0&1&0&0&0&0\\
	0&0&0&0&0&0&0&0&0\\
\end{array}
\right)\ .$$
\ee

Having given an intuitive explanation of our approach together with a worked example, we now formally describe the edges and their lengths. We include edges from the vertex
$u$ to the vertices $(1,t_1,d_1,1,r_1,S_1=t_1,R_1=r_1)$ for $t_1\in[n]$, $d_1\in[m]$,
$1\leq r_1\leq t_1$, of length $\sum\{g_j(0):t_1<j\leq t_0\}$.
Consider now any $h\in[k]$. For each $i$ with $d_{h-1}+1\leq i\leq d_h-1$ if any we include the
following edges of length $f_i(r_i)$ provided $1\leq r_{i+1}\leq r_i$ and $R_{i+1}\leq S_{i+1}$,
$$(h,t_h,d_h,i,r_i,S_i,R_i)\longrightarrow
(h,t_h,d_h,i+1,r_{i+1},S_{i+1}=S_i+t_h,R_{i+1}=R_i+r_{i+1})\ .$$
In addition, for every $h\in[k-1]$ and $i=d_h$ we include the following edges of length
$f_{d_h}(r_{d_h})+\sum\{g_j(d_h):t_{h+1}<j\leq t_h\}$
provided $1\leq r_{i+1}\leq r_i$ and $R_{i+1}\leq S_{i+1}$,
$$(h,t_h,d_h,i,r_i,S_i,R_i)\longrightarrow
(h+1,t_{h+1},d_{h+1},i+1,r_{i+1},S_{i+1}=S_i+t_{h+1},R_{i+1}=R_i+r_{i+1})\ .$$

Now consider $h=k$ and $i=d_h=d_k$. If $i=m$ then, provided $R_m=S_m$, we include the edges
$(k,t_k,m,m,r_m,S_m,R_m)\rightarrow v$ of length
$f_{d_h}(r_{d_h})+\sum\{g_j(m):0=t_{k+1}<j\leq t_k\}$.

Suppose now $h=k$ but $i=d_h=d_k<m$. We include the following edges
of length $f_{d_k}(r_{d_k})+\sum\{g_j(d_k):0=t_{k+1}<j\leq t_k\}$
provided $1\leq r_{i+1}\leq r_i$ and $R_{i+1}\leq S_{i+1}$,
$$(k,t_k,d_k,i,r_i,S_i,R_i)\longrightarrow
(k+1,t_{k+1}=0,d_{k+1}=m,i+1,r_{i+1},S_{i+1}=S_i,R_{i+1}=R_i+r_{i+1})\ .$$
Also, for each $i$ with $d_k+1\leq i\leq d_k-1=m-1$ if any we include the
following edges of length $f_i(r_i)$ provided $1\leq r_{i+1}\leq r_i$ and $R_{i+1}\leq S_{i+1}$,
$$(k+1,0,m,i,r_i,S_i,R_i)\longrightarrow
(k+1,0,m,i+1,r_{i+1},S_{i+1}=S_i,R_{i+1}=R_i+r_{i+1})\ .$$
Finally, if $R_m=S_m$, we include the edges
$(k+1,0,m,m,r_m,S_m,R_m)\rightarrow v$ of length $f_m(r_m)$.

\vskip.2cm
Now, as explained above, it is clear that each $r,c$ with $c$ of type $k$ and $r$
majorized by the conjugate $s$ of $c$ give a $u-v$ path of length $f(r)+g(c)$ in $D_k$.
Conversely, it is clear that every $u-v$ path in $D_k$ visits $m$ intermediate vertices
with $i=1,\dots,m$ and we can read off from this path $r=(r_1,\dots,r_m)$ directly and
$c=(c_1,\dots,c_n)$ and its conjugate $s=(s_1,\dots,s_m)$ of type $k$ as in \eqref{type}
with $r$ majorized by $s$, and $f(r)+g(c)$ equals the length of the path. So a shortest
directed $u-v$ path in $D_k$ gives a pair $r,c$ with $c$ of type $k$ minimizing $f(r)+g(c)$.

\vskip.2cm
Now, the number of vertices of $D_k$ is $O(kn^4m^4)$ and hence is polynomial in $m,n$. So a
shortest directed $u-v$ path in $D_k$ can be obtained in polynomial time, see e.g. \cite{Sch}.

\vskip.2cm
Now for $k=1,\dots,\min\{m,n\}$ we find the shortest path in $D_k$, read off $r,c$ with
minimum $f(r)+g(c)$ among those with $c$ of type $k$, compare to $c=(0^n)$ of type $k=0$
and $r=(0^m)$, and let $r,c$ be the best over all. We now use Lemma \ref{type} to obtain
a monotone matrix $A\in\{0,1\}^{m\times n}$ which has row and column sums $r,c$,
which is an optimal solution to our problem.
\epr

\section*{Acknowledgments}
The first author was partially supported by Charles University project
UNCE/SCI/004 and by the project 19-27871X of GA \v{C}R.
The second author was partially supported by a grant from the
Israel Science Foundation and by the Dresner chair at the Technion.

\end{document}